\documentclass[11pt,a4paper,leqno]{amsart}
\usepackage{lmodern}
\usepackage{amsmath,amssymb,amsthm}
\usepackage[expansion=false]{microtype}
\usepackage{booktabs}
\usepackage{hyperref}

\calclayout
\allowdisplaybreaks[1]
\newtheorem{theorem}{Theorem}[section]
\newtheorem{lemma}[theorem]{Lemma}
\newtheorem{proposition}[theorem]{Proposition}
\newtheorem{corollary}[theorem]{Corollary}
\theoremstyle{definition}

\theoremstyle{remark}
\newtheorem{remark}[theorem]{Remark}
\numberwithin{equation}{section}

\newcommand{\End}{\operatorname{End}}
\newcommand{\Hom}{\operatorname{Hom}}
\newcommand{\Ann}{\operatorname{Ann}}
\newcommand{\im}{\operatorname{im}}
\newcommand{\rank}{\operatorname{rank}}
\newcommand{\mrk}{\operatorname{mrk}}
\newcommand{\Span}{\operatorname{span}}
\newcommand{\vol}{\operatorname{vol}}
\newcommand{\ini}{\operatorname{in}}
\newcommand{\gr}{\operatorname{gr}}
\newcommand{\RR}{\mathbb R}

\newcommand{\cceil}[1]{\left\lceil #1\right\rceil}

\title[Linear bounds for commuting matrices]{Linear bound for commuting matrix algebra II}

    \author{Caihua Luo}
    \address{SCHOOL OF SCIENCE AND ENGINEERING, THE CHINESE UNIVERSITY OF HONG KONG, SHENZHEN,
GUANGDONG 518172, CHINA.}
    \email{luocaihua@cuhk.edu.cn}

    \author{Xu Song}
    \address{SCHOOL OF MATHEMATICAL SCIENCES, HEBEI NORMAL UNIVERSITY, SHIJIAZHUANG, HEBEI 050024, CHINA.}
    \email{songx@hebtu.edu.cn}

    \author{Yilin Zeng}
    \address{SCHOOL OF SCIENCE AND ENGINEERING, THE CHINESE UNIVERSITY OF HONG KONG, SHENZHEN,
GUANGDONG 518172, CHINA.}
    \email{223015132@link.cuhk.edu.cn}

\hypersetup{pdftitle={Linear Bounds for Commuting Matrix Algebras II},
pdfsubject={Dimension bounds for algebras generated by commuting matrices},
pdfauthor={},pdfkeywords={commuting matrices, regularized rank, faithful modules}}

\begin{document}
\begin{abstract}
Let $\mathbb{F}$ be a field, and let $A_1,\cdots,A_m$ be pairwise commuting
$n\times n$ matrices. For $m\geq3$ and $n\geq2$, we prove
\[
 \dim_{\mathbb{F}} \mathbb{F}[A_1,\cdots,A_m]
 \leq \left\lceil C_m(n-1)\right\rceil \quad\text{with}\quad
 C_m=\frac{\Gamma(m)}{\Gamma((m+1)/2)^2}.
\]
In particular, three
commuting matrices generate an algebra of dimension at most $2n-2$.
\end{abstract}

\maketitle
\markboth{LINEAR BOUNDS FOR COMMUTING MATRICES}{LINEAR BOUNDS FOR COMMUTING MATRICES}

\section{Introduction}

Let $\mathbb{F}$ be a field. For pairwise commuting matrices
$A_1,\cdots,A_m\in M_n(\mathbb{F})$, we consider the dimension of the unital algebra $\mathbb{F}[A_1,\cdots,A_m]$.

For one generator, the Cayley--Hamilton theorem gives the bound $n$. For two generators, the same bound follows from the commuting-pair theorem of Gerstenhaber \cite{Gerstenhaber}; see also the module-theoretic proof of Wadsworth \cite{Wadsworth}. The corresponding assertion for
three generators is the Gerstenhaber problem. For a fixed number of
generators, Bergman asked whether a linear bound exists
\cite[\S6, Question 15]{Bergman}. Li, Luo, Song, and Zeng
\cite{LLSZ} proved
\begin{equation}\label{eq:oldbound}
 \dim_{\mathbb{F}} \mathbb{F}[A_1,\cdots,A_m]\leq(2^m-1)(n-1)+1,
\end{equation}
and determined the exponential growth rate of the optimal coefficients.

In this paper, we improve the linear coefficient in Formula
\eqref{eq:oldbound} as follows. For $m\geq2$, set, $\Gamma$ denotes the classical gamma function,
\begin{equation}\label{eq:Cm}
 C_m:=\frac{\Gamma(m)}{\Gamma((m+1)/2)^2}.
\end{equation}

\begin{theorem}\label{thm:main}
Let $\mathbb{F}$ be a field, let $m\geq3$ and $n\geq2$, and let
$A_1,\cdots,A_m\in M_n(\mathbb{F})$ commute pairwise. Then
\begin{equation}\label{eq:mainstrict}
 \dim_{\mathbb{F}} \mathbb{F}[A_1,\cdots,A_m]<C_m(n-1)+1.
\end{equation}
Consequently,
\begin{equation}\label{eq:mainceil}
 \dim_{\mathbb{F}} \mathbb{F}[A_1,\cdots,A_m]\leq\cceil{C_m(n-1)}.
\end{equation}
In particular,
\[
 \dim_{\mathbb{F}} \mathbb{F}[A,B,C]\leq2n-2,
 \qquad
 \dim_{\mathbb{F}} \mathbb{F}[A_1,\cdots,A_5]\leq6n-6.
\]
\end{theorem}

\subsection{Multiplication pairings and regularized rank}
For a nonzero finite-dimensional commutative $\mathbb{F}$-algebra $R$, set
$R^*=\Hom_{\mathbb{F}}(R,\mathbb{F})$ and define
\begin{equation}\label{eq:pairings}
 \beta_\lambda(x,y)=\lambda(xy),\qquad
 g_{\mathbb{F}}(R)=\max_{\lambda\in R^*}\rank\beta_\lambda.
\end{equation}
We regard $\beta_\lambda$ also as a map $R\to R^*$ and write
$S_R=\{\beta_\lambda:\lambda\in R^*\}$. Its regularized maximum rank is defined as
\begin{equation}\label{eq:regg}
 \widetilde{g}_{\mathbb{F}}(R)=\lim_{q\to\infty}g_{\mathbb{F}}(R^{\otimes q})^{1/q} \quad \text{(limit exists by supermultiplicativity)}.
\end{equation}
 The next theorem relates this invariant to faithful-module dimension and identifies it with the numerical data of a canonical filtration. We construct the filtration, prove its uniqueness and duality, and establish the required tensor-product estimates in Sections~\ref{sec:expansion}--\ref{sec:exactrank}.

\begin{theorem}\label{thm:rank}
Let $R$ be a nonzero finite-dimensional commutative unital $\mathbb{F}$-algebra.
If $R$ acts faithfully on an $n$-dimensional space, then
\begin{equation}\label{eq:faithfulregg}
 g_{\mathbb{F}}(R)\leq\widetilde{g}_{\mathbb{F}}(R)\leq n.
\end{equation}
If $\mathbb{F}$ is infinite and $(a_i,b_i)_{i=1}^s$ are the dimensions of the codomain and domain of the canonical blocks of $S_R$, then
\begin{equation}\label{eq:exactregg}
 \widetilde{g}_{\mathbb{F}}(R)=\sum_{i=1}^s\sqrt{a_i b_i}.
\end{equation}
If $\mathbb{F}$ is infinite, $R$ is generated by at most $m\geq3$ elements,
and $D=\dim_{\mathbb{F}} R>1$, then
\begin{equation}\label{eq:affineregg}
 \widetilde{g}_{\mathbb{F}}(R)>1+\frac{D-1}{C_m}.
\end{equation}
\end{theorem}

\begin{theorem}\label{thm:optimal}
Fix an infinite field $\mathbb{F}$ and an integer $m\geq2$. Then
\[
 \sup_R\frac{\dim_{\mathbb{F}} R}{\widetilde{g}_{\mathbb{F}}(R)}=C_m,
\]
where $R$ ranges over the nonzero finite-dimensional commutative unital
$\mathbb{F}$-algebras generated by at most $m$ elements.
\end{theorem}

\subsection{The generator constants}
Define
\begin{equation}\label{eq:cm}
 c_m(\mathbb{F})=\sup\left\{
 \frac{\dim_{\mathbb{F}} \mathbb{F}[A_1,\cdots,A_m]}{n}:
 \begin{array}{c}
 n\geq1,\ A_i\in M_n(\mathbb{F}),\\[-2pt]
 A_iA_j=A_jA_i\text{ for all }i,j
 \end{array}\right\}.
\end{equation}
The single-matrix and commuting-pair results give $c_1(\mathbb{F})=c_2(\mathbb{F})=1$.
For $m\geq3$, Theorem~\ref{thm:main} gives $c_m(\mathbb{F})\leq C_m$.

\begin{corollary}\label{cor:growth}
For every field $\mathbb{F}$ and $m\geq4$,
\begin{equation}\label{eq:twosided}
 L_m:=\binom{m-2}{\lfloor(m-2)/2\rfloor}
 \leq c_m(\mathbb{F})\leq C_m.
\end{equation}
Furthermore,
\[
 C_m\sim\frac{2^m}{\sqrt{2\pi m}},\qquad
 \frac{C_m}{L_m}\longrightarrow2.
\]
\end{corollary}

The lower bound in Formula \eqref{eq:twosided} is the construction of
\cite[\S3.3]{LLSZ}, and the upper bound follows from
Theorem~\ref{thm:main}. Neither
Theorem~\ref{thm:optimal} nor Corollary~\ref{cor:growth} determines
$c_m(\mathbb{F})$ for $m\geq3$.

\subsection{Outline of the proof}
The proof begins with a maximal cyclic subspace. Its annihilator has square zero, giving a preliminary pairing-rank bound with a universal factor of two. Faithful tensor powers remove that factor and yield $\widetilde{g}_{\mathbb{F}}(R)\leq n$.

We then construct the canonical slope filtration of a matrix space. Tensor products preserve its rectangular expansion inequalities. For a symmetric space, duality pairs the block dimensions; tensor words with symmetric frequencies give an asymptotic rank lower bound $\sum_i\sqrt{a_i b_i}$. A complementary cover of the rows and columns gives the matching upper bound.

Finally, the product-colength inequality of \cite{LLSZ} places the normalized filtration polygon above $F_m(u)=[1-(1-u)^{1/m}]^m$. The integral of $\sqrt{F_m'}$ is $1/C_m$. The positive integral block dimensions bound both extreme slopes of the polygon. Retaining the resulting contributions at both ends of a truncated-derivative integral gives Formula \eqref{eq:affineregg}. Scalar extension to $\mathbb{F}(t)$ proves the matrix bound over finite fields. Truncated polynomial algebras then establish Theorem~\ref{thm:optimal}.

Throughout, all vector spaces are finite-dimensional, all algebras and modules are unital, and all algebras under consideration are commutative. Tensor products and duals are over the ground field. Field subscripts are omitted when unambiguous, except in $\widetilde{g}_{\mathbb{F}}(R)$. The notation $\perp$ always denotes an algebraic annihilator under the natural dual pairing.

\section{Faithful modules and multiplication pairings}\label{sec:faithful}

\subsection{Maximal rank and maximal cyclic subspaces}
We first record the elementary rank-pencil argument used below. It is a form of the Flanders--Atkinson lemma (see \cite[Lemma 2.2]{Pazzis}).

\begin{lemma}\label{lem:pencil}
Suppose $\mathbb{F}$ is infinite and $M$ has maximum rank in a vector space
$\mathcal T\subseteq\Hom_{\mathbb{F}}(U,W)$. Then
\[
 T(\ker M)\subseteq\im M\qquad(T\in\mathcal T).
\]
\end{lemma}

\begin{proof}
Put $r=\rank M$ and choose bases in which $M$ consists of an
$r\times r$ identity block and zeros elsewhere. Every
$(r+1)\times(r+1)$ minor of $M+tT$ vanishes for every $t\in \mathbb{F}$,
and therefore vanishes identically. The minor using the identity block
and one additional row and column has linear coefficient equal to the
corresponding entry of the lower-right block of $T$. That block is
zero, as required. If $M$ is injective or surjective, the assertion is
immediate.
\end{proof}

\begin{lemma}\label{lem:cyclic}
Suppose $\mathbb{F}$ is infinite and $R\subseteq\End_{\mathbb{F}}(V)$. Choose $v\in V$
so that $r=\dim_{\mathbb{F}} Rv$ is maximal, and put $J=\Ann_R(v)$. Then
\begin{equation}\label{eq:cyclic}
 JV\subseteq Rv,\qquad J(Rv)=0,\qquad J^2=0,
 \qquad\dim_{\mathbb{F}} R/J=r.
\end{equation}
\end{lemma}

\begin{proof}
Consider the vector space of evaluation maps
$\varepsilon_w:R\to V$, $a\mapsto aw$, as $w$ varies through $V$.
The map $\varepsilon_v$ has maximum rank. By
Lemma~\ref{lem:pencil}, $\varepsilon_w(J)\subseteq Rv$ for every $w$,
which gives $JV\subseteq Rv$. Commutativity gives $J(Rv)=0$.
It follows that $J^2V=0$, and faithfulness gives $J^2=0$.
The evaluation map at $v$ identifies $R/J$ with $Rv$.
\end{proof}

\begin{proposition}\label{prop:preliminary}
If $\mathbb{F}$ is infinite and $R$ acts faithfully on an $n$-dimensional
space, then $g_{\mathbb{F}}(R)\leq2n$.
\end{proposition}

\begin{proof}
Identify $R$ with its image in $\End_{\mathbb{F}}(V)$ and use the ideal $J$
from Lemma~\ref{lem:cyclic}. Choose a vector-space complement $E$
to $J$ in $R$, so $\dim E=r\leq n$. Since $J^2=0$, every
multiplication pairing has a matrix of the form
\[
 \begin{pmatrix} B&C\\ C^{\mathsf T}&0\end{pmatrix}
\]
relative to $R=E\oplus J$. The first $r$ columns contribute rank
at most $r$, and all remaining columns are supported on the first
$r$ rows. Hence the total rank is at most $2r\leq 2n$.
\end{proof}

\subsection{Tensor powers}
The preceding constant is independent of the number of generators,
which permits its removal by tensoring.

\begin{lemma}\label{lem:faithfultensor}
If $R\subseteq\End_{\mathbb{F}}(V)$ and $T\subseteq\End_{\mathbb{F}}(W)$, then the
natural action of $R\otimes T$ on $V\otimes W$ is faithful.
Moreover,
\begin{equation}\label{eq:tensorpairing}
 S_{R\otimes T}=S_R\otimes S_T,
 \qquad g_{\mathbb{F}}(R\otimes T)\geq g_{\mathbb{F}}(R)g_{\mathbb{F}}(T).
\end{equation}
\end{lemma}

\begin{proof}
Tensoring the two injective linear maps into the endomorphism spaces
remains injective over a field. As tensor products of matrix units are the matrix
units for the tensor-product basis, so the natural map
\[
 \End_{\mathbb{F}}(V)\otimes\End_{\mathbb{F}}(W)\longrightarrow\End_{\mathbb{F}}(V\otimes W)
\]
is an isomorphism. This proves faithfulness.

For $\lambda\in R^*$ and $\mu\in T^*$, direct evaluation on pure
tensors gives $\beta_{\lambda\otimes\mu}=\beta_\lambda\otimes
\beta_\mu$. Note that pure tensor functionals span $(R\otimes T)^*$, so
\[
 S_{R\otimes T}=S_R\otimes S_T.
\]
As the rank of a tensor product of matrices is the product of the ranks, so
\[
 g_{\mathbb{F}}(R\otimes T)\geq g_{\mathbb{F}}(R)g_{\mathbb{F}}(T).
\]
\end{proof}

\begin{proposition}\label{prop:faithfulregg}
For every field $\mathbb{F}$, if $R$ acts faithfully on an $n$-dimensional
space, then the limit in Formula \eqref{eq:regg} exists and
$g_{\mathbb{F}}(R)\leq\widetilde{g}_{\mathbb{F}}(R)\leq n$.
\end{proposition}

\begin{proof}
First, let $\mathbb{F}$ be infinite. Fix $\lambda\in R^*$ and write
$h=\rank\beta_\lambda$. As the algebra $R^{\otimes q}$ acts faithfully on $V^{\otimes q}$, and the pairing attached to $\lambda^{\otimes q}$ has rank $h^q$. Then Proposition~\ref{prop:preliminary} gives $h^q\leq2n^q$ for every $q\geq1$. Taking $q$-th roots and letting $q\to\infty$ yields $h\leq n$.

For an arbitrary field, extend scalars to the infinite field $K=\mathbb{F}(t)$. The extended action remains faithful and has dimension $n$ over $K$. Each original pairing has the same matrix after extension, so its rank is unchanged. The infinite-field argument therefore gives $g_{\mathbb{F}}(R)\leq n$ over every field.

Put $g_q=g_{\mathbb{F}}(R^{\otimes q})$ and $D=\dim_{\mathbb{F}} R$. By Lemma~\ref{lem:faithfultensor}, $g_{p+q}\geq g_pg_q$, and $1\leq g_q\leq D^q$. The superadditive-sequence lemma applied to $\log g_q$ gives
\[
 \lim_{q\to\infty}g_q^{1/q}=\sup_{q\geq1}g_q^{1/q}.
\]
In particular, this limit is at least $g_{\mathbb{F}}(R)$. Applying the already
proved bound $g\leq\dim V$ to the faithful action on $V^{\otimes q}$
gives $g_q\leq n^q$, and hence $\widetilde{g}_{\mathbb{F}}(R)\leq n$.
\end{proof}

\section{Expansion of matrix spaces}\label{sec:expansion}

For a subspace $S\subseteq\Hom_{\mathbb{F}}(U,W)$ and a subspace $X\subseteq U$,
write
\[
 S(X)=\Span_{\mathbb{F}}\{s(x):s\in S,\ x\in X\}.
\]
If $U\ne0$, define the minimum expansion
\[
 e(S)=\min_{0\ne X\subseteq U}\frac{\dim S(X)}{\dim X}.
\]
The minimum exists because there are only finitely many possible pairs
of dimensions. We say that $S:U\to W$ is \emph{ratio-balanced} if
$b=\dim U>0$, $a=\dim W>0$, and
\begin{equation}\label{eq:balanced}
 \dim S(X)\geq\frac ab\dim X\qquad(X\subseteq U).
\end{equation}
Then $S(U)=W$, the common kernel is zero, and $e(S)=a/b$.

\begin{lemma}\label{lem:amplification}
For every nonzero finite-dimensional space $E$,
$e(S\otimes I_E)=e(S)$.
\end{lemma}

\begin{proof}
Choose a complete flag $0=E_0<E_1<\cdots<E_d=E$. For an arbitrary
subspace $X\subseteq U\otimes E$, set
\[
 X_j=X\cap(U\otimes E_j),\qquad
 Y=(S\otimes I_E)(X),\qquad Y_j=Y\cap(W\otimes E_j).
\]
After choosing $E_j/E_{j-1}\cong \mathbb{F}$, the successive quotients of these filtrations identify with subspaces $\overline X_j\subseteq U$ and $\overline Y_j\subseteq W$. Since every $s\otimes I_E$ preserves the flags, $S(\overline X_j)\subseteq\overline Y_j$. Thus
\[
 \dim Y\geq\sum_j\dim S(\overline X_j)
 \geq e(S)\sum_j\dim\overline X_j=e(S)\dim X.
\]
The reverse inequality follows from taking a minimizing subspace $X_0\subseteq U$ and a nonzero $v\in E$, and computing
\[
\frac{\dim S(X_0\otimes\mathbb{F}v)}{\dim(X_0\otimes\mathbb{F}v)}.
\]
\end{proof}

\begin{proposition}\label{prop:tensorexpansion}
Let $S: U\to W$ and $T: U'\to W'$ be matrix spaces with nonzero domains. Then
\begin{equation}\label{eq:tensorexpansion}
 e(S\otimes T)=e(S)e(T).
\end{equation}
In particular, tensor products of ratio-balanced spaces are ratio-balanced at the product dimension ratio.
\end{proposition}

\begin{proof}
If the codomain is zero, the assertion is immediate. Otherwise,
for $X\subseteq U\otimes U'$, put $Y=(S\otimes I_{U'})(X)$.
As linear spans,
\[
 (S\otimes T)(X)=(I_W\otimes T)(Y).
\]
Lemma~\ref{lem:amplification}, also applied after interchanging the
factors, gives
\[
 \dim(S\otimes T)(X)\geq e(T)\dim Y
 \geq e(T)e(S)\dim X.
\]
This remains true when $Y=0$. For minimizing subspaces $X_0\subseteq U$
and $X'_0\subseteq U'$, one has
$(S\otimes T)(X_0\otimes X'_0)=S(X_0)\otimes T(X'_0)$.
The corresponding ratio gives the reverse inequality.
\end{proof}

\begin{lemma}\label{lem:halfrank}
Suppose $\mathbb{F}$ is infinite. If $T\subseteq\Hom_{\mathbb{F}}(U,W)$ is
ratio-balanced and $\dim U=\dim W=h$, then $\mrk(T)\geq h/2$,
where $\mrk(T)=\max_{t\in T}\rank t$.
\end{lemma}

\begin{proof}
Choose $M\in T$ of maximum rank $r$. By Lemma~\ref{lem:pencil},
$T(\ker M)\subseteq\im M$. The expansion condition gives
\[
 h-r=\dim\ker M\leq\dim T(\ker M)\leq r.\qedhere
\]
\end{proof}

\section{The canonical slope filtration}\label{sec:filtration}

We say that $S\subseteq\Hom_{\mathbb{F}}(U,W)$ has \emph{full support} if
\[
 \bigcap_{s\in S}\ker s=0\quad\text{and}\quad S(U)=W.
\]
The filtration below is the generalized Kronecker instance of the Harder--Narasimhan filtration (see \cite[\S1.3]{IwamasaOkiSoma}). We give its construction and duality directly over $\mathbb{F}$ as follows.

\begin{proposition}\label{prop:filtration}
Let $S:U\to W$ have full support, with $U,W\ne0$. There are unique
filtrations
\begin{equation}\label{eq:flags}
 \begin{split}
 0=U_0<U_1<\cdots<U_s=U,\\
 0=W_0<W_1<\cdots<W_s=W
 \end{split}
\end{equation}
such that $W_i=S(U_i)$, the induced diagonal spaces
\[
 S_i\subseteq\Hom_{\mathbb{F}}(U_i/U_{i-1},W_i/W_{i-1})
\]
are ratio-balanced. Their dimensions
\[
 b_i=\dim(U_i/U_{i-1}),\qquad a_i=\dim(W_i/W_{i-1})
\]
satisfy
\begin{equation}\label{eq:slopes}
 a_i,b_i>0,\qquad
 \frac{a_1}{b_1}<\frac{a_2}{b_2}<\cdots<\frac{a_s}{b_s}.
\end{equation}
\end{proposition}

\begin{proof}
Let $\alpha=e(S)>0$. As the function $f(X)=\dim S(X)-\alpha\dim X$ is nonnegative on subspaces of $U$, so the relations
\[
 S(X+Y)=S(X)+S(Y),\qquad S(X\cap Y)\subseteq S(X)\cap S(Y)
\]
show that $f(X+Y)+f(X\cap Y)\leq f(X)+f(Y)$. Hence sums and intersections of minimizing subspaces remain minimizing. Finite dimensionality gives a unique largest minimizer $U_1\ne0$. Put $W_1=S(U_1)$, then the restriction to $U_1\to W_1$ is ratio-balanced with ratio $\alpha$.

Suppose $U_1\ne U$. The induced quotient space $\overline S$ maps from $U/U_1$ to $W/W_1$. If $X\ne0$ is a quotient subspace with full inverse image $\widetilde X$, then
\[
 \dim\widetilde X=\dim U_1+\dim X,\qquad
 \dim S(\widetilde X)=\dim W_1+\dim\overline S(X).
\]
Notice that a quotient expansion ratio smaller than $\alpha$ would contradict minimality, and equality would contradict maximality of $U_1$. Thus $e(\overline S)>\alpha$. In particular, the quotient again has full support. Iteration constructs Formula \eqref{eq:flags}--\eqref{eq:slopes}.

For uniqueness, start with any filtration having the stated
properties. Filter a subspace $X\subseteq U$ by $X\cap U_i$ and
its image by $S(X)\cap W_i$. The diagonal expansion inequalities
give, let $\gr_i(X) = \frac{X\cap U_i}{X\cap U_{i-1}}$,
\[
 \dim S(X)\geq\sum_i\frac{a_i}{b_i}\dim\gr_i X.
\]
The first slope is therefore $e(S)$. Equality at that slope forces
$\gr_iX=0$ for $i>1$, so every minimizer is contained in $U_1$.
Thus $U_1$ is the unique largest minimizer, and $W_1=S(U_1)$ is
determined. Apply the same argument to the quotient.
\end{proof}

For a subspace $T\subseteq\Hom_{\mathbb{F}}(U,W)$, write
\[
 T^*=\{t^*:t\in T\}\subseteq\Hom_{\mathbb{F}}(W^*,U^*)
\]
for its space of dual maps.

\begin{lemma}\label{lem:dualbalanced}
If $T:U\to W$ is ratio-balanced with dimensions $b,a>0$, then
$T^*:W^*\to U^*$ is ratio-balanced with ratio $b/a$.
\end{lemma}

\begin{proof}
For $X\subseteq W^*$, put $x=\dim X$ and $y=\dim T^*(X)$.
One has $T((T^*(X))^\perp)\subseteq X^\perp$. Applying Formula \eqref{eq:balanced} to the subspace $(T^*(X))^\perp$ of dimension $b-y$ yields
$(a/b)(b-y)\leq a-x$, i.e., $y\geq(b/a)x$.
\end{proof}

\begin{proposition}\label{prop:symmetry}
The canonical filtration of $S^*:W^*\to U^*$ is
\[
 \widehat U_i=W_{s-i}^\perp,\qquad
 \widehat W_i=U_{s-i}^\perp\qquad(0\leq i\leq s).
\]
If $W=U^*$ and every member of $S$ is symmetric, then
\begin{equation}\label{eq:symmetry}
 U_i=W_{s-i}^\perp,\qquad W_i=U_{s-i}^\perp,
 \qquad a_i=b_{s+1-i},\quad b_i=a_{s+1-i}.
\end{equation}
\end{proposition}

\begin{proof}
Compatibility follows by dualizing $S(U_j)\subseteq W_j$. The diagonal spaces are the duals of the original blocks in reverse order. Lemma~\ref{lem:dualbalanced} makes them ratio-balanced, with strictly increasing slopes $b_s/a_s<\cdots<b_1/a_1$. As the preceding flag term $\widehat{W}_{i-1}$ is already in the image, and the image on the current quotient is full, so their full quotient images show inductively that $S^*(\widehat U_i)=\widehat W_i$. Proposition~\ref{prop:filtration} identifies these as the canonical flags. In the symmetric case, the natural identification $U^{**}=U$ identifies $S^*$ with $S$. Uniqueness gives the flag identities, and taking dimension differences proves the remaining assertions.
\end{proof}

\section{An exact formula for regularized rank}\label{sec:exactrank}

For a nonzero matrix space $S$, put
\[
 \rho_{\mathbb{F}}(S)=\lim_{q\to\infty}\mrk(S^{\otimes q})^{1/q}.
\]
The limit exists by the same supermultiplicative argument as in Proposition~\ref{prop:faithfulregg}. On the other hand, Lemma~\ref{lem:faithfultensor} gives $\rho_{\mathbb{F}}(S_R)=\widetilde{g}_{\mathbb{F}}(R)$.

\subsection{An upper bound from a cover of rows and columns}

\begin{proposition}\label{prop:cover}
Suppose every member of $S$ is block upper triangular, with positive
row-block dimensions $a_1,\cdots,a_s$, positive column-block
dimensions $b_1,\cdots,b_s$, and nondecreasing ratios $a_i/b_i$.
Then, over any field and for every $q\geq1$,
\begin{equation}\label{eq:cover}
 \mrk(S^{\otimes q})\leq
 \left(\sum_{i=1}^s\sqrt{a_i b_i}\right)^q.
\end{equation}
\end{proposition}

\begin{proof}
Set $w_i=\sqrt{a_i/b_i}$. For a word
$\boldsymbol i=(i_1,\cdots,i_q)\in\{1,\cdots,s\}^q$, define
\[
 A_{\boldsymbol i}=\prod_\nu a_{i_\nu},\qquad
 B_{\boldsymbol i}=\prod_\nu b_{i_\nu},\qquad
 w_{\boldsymbol i}=\prod_\nu w_{i_\nu}.
\]
A nonzero tensor block from a column word $\boldsymbol j$ to a row
word $\boldsymbol i$ requires $i_\nu\leq j_\nu$ for every $\nu$.
Thus $w_{\boldsymbol i}\leq w_{\boldsymbol j}$.

Select all rows in blocks with $w_{\boldsymbol i}\leq1$ and all
columns in blocks with $w_{\boldsymbol j}>1$. Every nonzero entry
is covered. A matrix supported on a union of selected rows and
selected columns has rank at most their total number, by separating
the selected rows from the remaining rows. Moreover,
\[
 A_{\boldsymbol i}\leq\sqrt{A_{\boldsymbol i}B_{\boldsymbol i}}
 \quad(w_{\boldsymbol i}\leq1),\qquad
 B_{\boldsymbol i}\leq\sqrt{A_{\boldsymbol i}B_{\boldsymbol i}}
 \quad(w_{\boldsymbol i}>1).
\]
The two conditions partition the word labels. The total number of
selected rows and columns is therefore at most
\[
 \sum_{\boldsymbol i}\sqrt{A_{\boldsymbol i}B_{\boldsymbol i}}
 =\left(\sum_i\sqrt{a_i b_i}\right)^q.\qedhere
\]
\end{proof}

\subsection{The matching lower bound for symmetric spaces}

\begin{theorem}\label{thm:exactspace}
Let $\mathbb{F}$ be infinite, let $U\ne0$, and let
$S\subseteq\Hom_{\mathbb{F}}(U,U^*)$ have full support and consist of
symmetric maps. If $(a_i,b_i)_{i=1}^s$ are
its canonical block dimensions, then
\begin{equation}\label{eq:rhoformula}
 \rho_{\mathbb{F}}(S)=\sum_{i=1}^s\sqrt{a_i b_i}.
\end{equation}
\end{theorem}

\begin{proof}
Choose bases adapted to the canonical flags. Every member of $S$
is block upper triangular. Proposition~\ref{prop:cover} proves the
upper bound. For the lower bound, set
$c_i=\sqrt{a_i b_i}$ and $E=\sum_i c_i$.

Order the block words of $S^{\otimes q}$ lexicographically.
The resulting matrices are block upper triangular, and the diagonal
word $(i_1,\cdots,i_q)$ ranges over
$S_{i_1}\otimes\cdots\otimes S_{i_q}$.
Choose nonnegative integers $q_1,\cdots,q_s$ such that
\begin{equation}\label{eq:counts}
 \sum_iq_i=q,\qquad q_i=q_{s+1-i}.
\end{equation}
For a word having these multiplicities,
Proposition~\ref{prop:symmetry} makes the diagonal block square,
of size
\begin{equation}\label{eq:blocksize}
 h=\prod_i a_i^{q_i}=\prod_i b_i^{q_i}=\prod_i c_i^{q_i}.
\end{equation}
Proposition~\ref{prop:tensorexpansion} makes its matrix space
ratio-balanced with ratio one. Lemma~\ref{lem:halfrank} supplies
a minor of size at least $h/2$ that is not identically zero.

There are $q!/\prod_i q_i!$ words of this type. The corresponding
rank lower bounds can be attained simultaneously, even though
the blocks share coefficients. Indeed, every block projection is
the tensor product of the surjective maps $S\to S_i$, and is
therefore onto its diagonal block space. Pull back a nonzero minor
polynomial from each selected block to $S^{\otimes q}$. Their
product is nonzero and, since $\mathbb{F}$ is infinite, has a nonzero value
at some member of $S^{\otimes q}$. The selected rows and columns
form a square block upper triangular minor with nonsingular
diagonal blocks. Thus their ranks add, giving
\begin{equation}\label{eq:typelower}
 \mrk(S^{\otimes q})\geq
 \frac12\frac{q!}{\prod_iq_i!}\prod_i c_i^{q_i}.
\end{equation}

Put $p_i=c_i/E$. Symmetry gives $p_i=p_{s+1-i}$. Along even $q$,
choose symmetric integer counts as in Formula \eqref{eq:counts} with $q_i/q\to p_i$. Stirling's formula gives
\[
 \left(\frac{q!}{\prod_iq_i!}\right)^{1/q}
 \longrightarrow\exp\left(-\sum_i p_i\log p_i\right).
\]
The $q$-th root of the right side of Formula \eqref{eq:typelower} therefore tends to
\[
 \exp\left(-\sum_i p_i\log p_i\right)\prod_i c_i^{p_i}=E.
\]
Since the full limit defining $\rho_{\mathbb{F}}(S)$ exists, this subsequence
proves the lower bound.
\end{proof}

\begin{corollary}\label{cor:algebraexact}
Over an infinite field, every nonzero finite-dimensional commutative
unital algebra satisfies
$\widetilde{g}_{\mathbb{F}}(R)=\sum_i\sqrt{a_i b_i}$ for the canonical blocks of $S_R$.
\end{corollary}

\begin{proof}
Commutativity makes $S_R$ symmetric. The identities
$\beta_\lambda(1)=\lambda$ and $\beta_\lambda(x)(1)=\lambda(x)$
give full total image and zero common kernel. Apply
Theorem~\ref{thm:exactspace} and $\rho_{\mathbb{F}}(S_R)=\widetilde{g}_{\mathbb{F}}(R)$, we get $\widetilde{g}_{\mathbb{F}}(R)=\sum_i\sqrt{a_i b_i}$.
\end{proof}

\begin{remark}
Asymptotic maximum rank and its relation to shrinking subspaces are also studied in \cite{LiQiaoWangDuan}. The argument above uses all rectangular blocks and obtains simultaneous rank attainment by a single nonvanishing polynomial. It does not require independent coefficients on different diagonal blocks.    
\end{remark}

\section{Colength and the filtration polygon}\label{sec:polygon}

\subsection{The product-colength inequality}
For a finite-colength ideal $H$ in $P=\mathbb{F}[x_1,\cdots,x_m]$, put $\ell(H)=\dim_{\mathbb{F}}(P/H)$. The following product-colength inequality is due to Li, Luo, Song, and Zeng \cite{LLSZ}. We include its proof.

\begin{lemma}\label{lem:colength}
If $H,K\subseteq P$ have finite colength, then
\begin{equation}\label{eq:colength}
 \ell(HK)^{1/m}\leq\ell(H)^{1/m}+\ell(K)^{1/m}.
\end{equation}
\end{lemma}

\begin{proof}
If one ideal is $P$, the assertion is immediate. For $m=1$, it is
equality, since colength is the degree of a generator. Suppose
$m\geq2$ and first assume both ideals are monomial.
Let $C=\RR_{\geq0}^m$ and define
\[
 U_H=\bigcup_{x^\alpha\in H}(\alpha+C).
\]
The union can be taken over finitely many minimal monomial generators. Thus $U_H$ is closed and $U_H+C=U_H$. Notice that $y\in U_H$ is equivalent to $x^{\lfloor y\rfloor}\in H$. The complement is consequently the disjoint union of the half-open
unit cubes indexed by the standard monomials, so
$\vol(C\setminus U_H)=\ell(H)$.
Monomial multiplication gives $U_{HK}=U_H+U_K$.

Let $\Theta=S^{m-1}\cap C$, with surface measure $d\theta$.
On each ray, write
$U_H\cap\RR_{\geq0}\theta=[r_H(\theta),\infty)\theta$.
Finite colength of a monomial ideal implies that it contains a
power of each variable, so $r_H$ is finite and bounded. It is
measurable because
$\{\theta:r_H(\theta)\leq t\}=\{\theta:t\theta\in U_H\}$ is closed.
Polar integration yields
\[
 \ell(H)=\frac1m\int_\Theta r_H(\theta)^m\,d\theta.
\]
The Minkowski-sum identity implies $r_{HK}\leq r_H+r_K$. Minkowski's inequality in $L^m(\Theta,d\theta/m)$ proves Formula \eqref{eq:colength} for monomial ideals.

For general ideals, fix a monomial order. Standard monomials
give $\ell(H)=\ell(\ini H)$ and $\ell(K)=\ell(\ini K)$ (see \cite{CoxLittleOShea}). Leading monomials multiply, so
$(\ini H)(\ini K)\subseteq\ini(HK)$. The monomial product has
finite colength. Hence $HK$ also has finite colength and
\[
 \ell(HK)\leq\ell((\ini H)(\ini K))
 \leq\bigl(\ell(H)^{1/m}+\ell(K)^{1/m}\bigr)^m.\qedhere
\]
\end{proof}

\begin{lemma}\label{lem:annihilating}
Suppose $R$ is generated by at most $m$ elements and $D=\dim_{\mathbb{F}} R$.
For a subspace $U\subseteq R$, put $W=S_R(U)\subseteq R^*$. Then
\begin{equation}\label{eq:annihilating}
 D^{1/m}\leq(D-\dim U)^{1/m}+(\dim W)^{1/m}.
\end{equation}
\end{lemma}

\begin{proof}
If $u\in U$ and $v\in W^\perp$, then
$\lambda(uv)=\beta_\lambda(u)(v)=0$ for every $\lambda\in R^*$.
Thus $UW^\perp=0$. The ideals $I=RU$ and $J=RW^\perp$ satisfy
$IJ=0$. Present $R=P/Q$ with $P=\mathbb{F}[x_1,\cdots,x_m]$, and let
$H,K$ be the inverse images of $I,J$. Then $HK\subseteq Q$, while
\[
 \ell(H)=\dim R/I\leq D-\dim U,\qquad
 \ell(K)=\dim R/J\leq\dim W.
\]
Therefore, the lemma follows from Lemma~\ref{lem:colength} and the surjective map $P/HK\to P/Q$.
\end{proof}

\subsection{The normalized polygon}
For the canonical blocks of $S_R$, set
\[
 A_i=\sum_{j\leq i}a_j,\qquad B_i=\sum_{j\leq i}b_j,
 \qquad A_0=B_0=0.
\]
Both final sums equal $D$. Let $G$ be the piecewise linear function joining $(B_i/D,A_i/D)$. Its slopes $a_i/b_i$ are positive and strictly increasing. Hence $G$ is convex, increasing, and absolutely continuous, with $G(0)=0$ and $G(1)=1$. Lemma~\ref{lem:annihilating} at the filtration vertices gives
\[
 1\leq(1-B_i/D)^{1/m}+(A_i/D)^{1/m}.
\]
It follows that the vertices lie above
\begin{equation}\label{eq:Fm}
 F_m(u)=[1-(1-u)^{1/m}]^m.
\end{equation}
For $m\geq2$,
\begin{equation}\label{eq:Fprime}
 F_m'(u)=\bigl((1-u)^{-1/m}-1\bigr)^{m-1}.
\end{equation}
This derivative is increasing and integrable, so $F_m$ is convex
and absolutely continuous. Each segment of $G$ lies above $F_m$,
and therefore
\begin{equation}\label{eq:polygonabove}
 G\geq F_m,\qquad
 \frac1D\sum_i\sqrt{a_i b_i}=\int_0^1\sqrt{G'(u)}\,du.
\end{equation}

\begin{lemma}\label{lem:convexcomparison}
Let $F,G:[0,1]\to[0,1]$ be convex, increasing, absolutely continuous
functions with common endpoints $0,1$. For such a function $H$,
define
\[
 T_H(t)=\int_0^1\min\{H'(u),t\}\,du\qquad(t\geq0).
\]
Then
\begin{equation}\label{eq:truncated}
 T_H(t)=\min_{0\leq x\leq1}\{H(x)+t(1-x)\}.
\end{equation}
If $G\geq F$, then $T_G(t)\geq T_F(t)$ for every $t\geq0$, and
\begin{equation}\label{eq:integralcomparison}
 \int_0^1\sqrt{G'}-\int_0^1\sqrt{F'}
 =\frac14\int_0^\infty(T_G(t)-T_F(t))t^{-3/2}\,dt\geq0.
\end{equation}
If also $G'\geq\alpha>0$ almost everywhere and $F'$ is continuous
near zero with $F'(0)=0$, then the last inequality is strict.
\end{lemma}

\begin{proof}
The derivative of $H$ is nondecreasing. Thus the minimum in
Formula \eqref{eq:truncated} occurs where that derivative passes through
$t$, with the endpoint cases included. Splitting the integral there
proves the identity. The inequality for the minima follows from
$G\geq F$.

For $z\geq0$, direct integration gives
\[
 \sqrt z=\frac14\int_0^\infty\min\{z,t\}\,t^{-3/2}\,dt.
\]
Tonelli's theorem applied separately to $F'$ and $G'$ proves Formula \eqref{eq:integralcomparison}; the integrals are finite by Cauchy--Schwarz and $\int_0^1H'=1$.
For $0<t<\alpha$, one has $T_G(t)=t$, whereas continuity of
$F'$ near zero makes $T_F(t)<t$. Integration gives strictness.
\end{proof}

The substitution $v=(1-u)^{1/m}$ gives
\begin{equation}\label{eq:beta}
 \begin{split}
 I_m:=\int_0^1\sqrt{F_m'(u)}\,du
 &=m\int_0^1[v(1-v)]^{(m-1)/2}\,dv\\
 &=\frac{\Gamma((m+1)/2)^2}{\Gamma(m)}=\frac1{C_m}.
 \end{split}
\end{equation}
Since the initial slope of $G$ is positive, Lemma~\ref{lem:convexcomparison}
and Corollary~\ref{cor:algebraexact} give
\begin{equation}\label{eq:homogeneous}
 \widetilde{g}_{\mathbb{F}}(R)>D/C_m\qquad(m\geq2)
\end{equation}
over an infinite field. The next section quantifies the gap in this
inequality by using the integrality of the block dimensions. This
additional information yields the affine estimate in
Theorem~\ref{thm:main}.

\section{Endpoint estimates and the linear bound}\label{sec:endpoints}

The integers $a_i,b_i$ are positive and each is at most $D$.
Consequently the polygon satisfies
\begin{equation}\label{eq:slopebounds}
 D^{-1}\leq G'(u)\leq D\qquad\text{almost everywhere}.
\end{equation}
Both inequalities contribute to the refinement.

\begin{proposition}\label{prop:endpoint}
Let $m\geq2$ and $D\geq1$. Suppose $G$ is convex, increasing,
and absolutely continuous, with endpoints $0,1$, with $G\geq F_m$,
and with the derivative bounds in Formula \eqref{eq:slopebounds}. Set
\[
 E=D\int_0^1\sqrt{G'(u)}\,du,
 \qquad
 h_m(s)=1-(1+s^{2/(m-1)})^{-(m-1)}.
\]
Then
\begin{equation}\label{eq:tailgain}
 E-DI_m\geq D\int_0^{D^{-1/2}}h_m(s)\,ds.
\end{equation}
If $m\geq3$ and $D>1$, then
\begin{equation}\label{eq:endpointaffine}
 E>1+(D-1)I_m.
\end{equation}
\end{proposition}

\begin{proof}
Put $r=m-1$. To evaluate $T_{F_m}(t)$, substitute
$y=(1-x)^{1/m}$ in Formula \eqref{eq:truncated}. We must minimize
$(1-y)^m+ty^m$ on $[0,1]$. For $t>0$, its minimum occurs at
$y=(1+t^{1/r})^{-1}$ and equals
\begin{equation}\label{eq:Tformula}
 T_{F_m}(t)=\frac{t}{(1+t^{1/r})^r}.
\end{equation}
In particular,
\begin{equation}\label{eq:reciprocity}
 T_{F_m}(t)=tT_{F_m}(1/t)\qquad(t>0).
\end{equation}
The slope bounds give $T_G(t)=t$ for $0\leq t\leq D^{-1}$ and
$T_G(t)=1$ for $t\geq D$. In Formula \eqref{eq:integralcomparison}, discard
only the nonnegative contribution from $D^{-1}<t<D$.
The substitution $t\mapsto1/t$, together with Formula \eqref{eq:reciprocity}, identifies the two remaining integrals.
Therefore
\[
 \frac ED-I_m\geq\frac12\int_0^{D^{-1}}
 \left[t-\frac{t}{(1+t^{1/r})^r}\right]t^{-3/2}\,dt.
\]
Substituting $t=s^2$ proves Formula \eqref{eq:tailgain}.

For $m\geq3$, put $p=2/(m-1)\leq1$. Direct differentiation gives,
for $s>0$,
\begin{equation}\label{eq:concavity}
 h_m''(s)=2s^{p-2}(1+s^p)^{-(m+1)}
 \bigl[(p-1)-3s^p\bigr]<0.
\end{equation}
Thus $h_m$ is strictly concave on $[0,1]$, with $h_m(0)=0$.
We also have
\begin{equation}\label{eq:hint}
 \int_0^1h_m(s)\,ds=1-I_m.
\end{equation}
Indeed, take the straight-line function $G(u)=u$ and $D=1$ in
the preceding calculation. There is then no discarded middle
interval, so equality holds throughout and gives Formula \eqref{eq:hint}.

Now suppose $D>1$ and set $a=D^{-1/2}<1$.
Strict concavity implies $h_m(au)>a h_m(u)$ for $0<u\leq1$.
Since $Da^2=1$,
\[
 D\int_0^a h_m(s)\,ds
 =Da\int_0^1h_m(au)\,du
 >Da^2\int_0^1h_m(u)\,du=1-I_m.
\]
Combining this with Formula \eqref{eq:tailgain} proves Formula \eqref{eq:endpointaffine}.
\end{proof}

\begin{proof}[Completion of the proof of Theorem~\ref{thm:rank}]
Proposition~\ref{prop:faithfulregg} proves Formula \eqref{eq:faithfulregg},
and Corollary~\ref{cor:algebraexact} proves Formula \eqref{eq:exactregg}.
For the canonical polygon, Formula \eqref{eq:polygonabove} identifies $E=\sum_i\sqrt{a_i b_i}=\widetilde{g}_{\mathbb{F}}(R)$. Apply Proposition~\ref{prop:endpoint} and $I_m=1/C_m$ to obtain Formula \eqref{eq:affineregg}.
\end{proof}

\begin{proof}[Proof of Theorem~\ref{thm:main}]
First suppose $\mathbb{F}$ is infinite. Let $R$ act faithfully on an
$n$-dimensional space and write $D=\dim_{\mathbb{F}} R$. If $D>1$, then
Theorem~\ref{thm:rank} gives
\[
 n\geq\widetilde{g}_{\mathbb{F}}(R)>1+\frac{D-1}{C_m}.
\]
Rearranging proves Formula \eqref{eq:mainstrict}. If $D=1$ and $n\geq2$,
the same inequality is immediate.

For an arbitrary field, extend scalars to $K=\mathbb{F}(t)$. The algebra
$R_K=K\otimes_{\mathbb{F}}R$ is generated by at most $m$ elements and acts
faithfully on $K\otimes_{\mathbb{F}}V$. Its dimension over $K$ is $D$, and
the module dimension is $n$. Applying the infinite-field result
proves Formula \eqref{eq:mainstrict} over $\mathbb{F}$. This argument requires no
invariance statement about the canonical filtration or $\widetilde{g}_{\mathbb{F}}$.

If an integer $D$ satisfies $D<x+1$, then $D\leq\lceil x\rceil$.
Taking $x=C_m(n-1)$ proves Formula \eqref{eq:mainceil}.
For commuting matrices, use $R=\mathbb{F}[A_1,\cdots,A_m]\subseteq M_n(\mathbb{F})$.
Finally, a unital subalgebra of $M_1(\mathbb{F})=\mathbb{F}$ has dimension one.
\end{proof}

\begin{corollary}\label{cor:quanttriple}
Let $\mathbb{F}$ be infinite, and let $R$ be a nonzero finite-dimensional
commutative unital $\mathbb{F}$-algebra generated by at most three elements.
Then, with $D=\dim_{\mathbb{F}} R$,
\begin{equation}\label{eq:triplecorrection}
 \widetilde{g}_{\mathbb{F}}(R)\geq\frac D2+\frac{\sqrt D}{\sqrt D+1}.
\end{equation}
For $D>1$, the right side is strictly larger than $(D+1)/2$.
\end{corollary}

\begin{proof}
Here $I_3=1/2$ and $h_3(s)=1-(1+s)^{-2}$. Since
\[
 \int_0^a h_3(s)\,ds=\frac{a^2}{1+a},
\]
the assertion follows from Formula \eqref{eq:tailgain} with $a=D^{-1/2}$.
\end{proof}

\section{Optimality of the regularized rank coefficient}\label{sec:optimality}

Fix an infinite field $\mathbb{F}$ and $m\geq2$. Let
\[
 R_d=\mathbb{F}[x_1,\cdots,x_m]/(x_1,\cdots,x_m)^{d+1},
\]
\[
 D_d=\binom{d+m}{m},\qquad h_j=\binom{j+m-1}{m-1}.
\]
The number $h_j$ is the dimension of the degree-$j$ homogeneous
component. We obtain the needed upper rank bound directly from
these degree spaces, without identifying them with the canonical
filtration.

\begin{lemma}\label{lem:truncated}
For every $d\geq0$,
\begin{equation}\label{eq:truncatedupper}
 \widetilde{g}_{\mathbb{F}}(R_d)\leq\sum_{j=0}^d\sqrt{h_jh_{d-j}}.
\end{equation}
\end{lemma}

\begin{proof}
Order the degree components of $R_d$ from $d$ down to $0$,
and their duals in $R_d^*$ from $0$ up to $d$.
Index both by $j=0,\cdots,d$. The block dimensions are then
$a_j=h_j$ and $b_j=h_{d-j}$.
An entry from degree $d-j$ in $R_d$ to the dual of degree $i$
vanishes for $i>j$, since the product has degree greater
than $d$. Thus every multiplication pairing is block upper
triangular. Since $h_j$ increases with $j$, the ratios
$h_j/h_{d-j}$ are nondecreasing. Proposition~\ref{prop:cover}
applies to every tensor power. Taking roots and limits proves Formula \eqref{eq:truncatedupper}.
\end{proof}

\begin{proof}[Proof of Theorem~\ref{thm:optimal}]
The Formula \eqref{eq:homogeneous} gives
$\dim_{\mathbb{F}} R/\widetilde{g}_{\mathbb{F}}(R)<C_m$ for $R$.
It remains to show that truncated polynomial algebras approach
this constant.

For fixed $m$, $h_j$ is a polynomial in $j$ of degree $m-1$ with
leading coefficient $1/(m-1)!$. Uniformly for $0\leq j\leq d$,
\[
 \frac{h_j}{d^{m-1}}-
 \frac{(j/d)^{m-1}}{(m-1)!}\longrightarrow0.
\]
Uniform continuity of the square root on bounded nonnegative
intervals then justifies the Riemann-sum limit
\begin{equation}\label{eq:riemann}
 \begin{split}
 \frac{\sum_{j=0}^d\sqrt{h_jh_{d-j}}}{D_d}
 &\longrightarrow
 m\int_0^1[t(1-t)]^{(m-1)/2}\,dt\\
 &=I_m=\frac1{C_m},
 \end{split}
\end{equation}
where we used $D_d\sim d^m/m!$.
Combining Formula \eqref{eq:homogeneous} and Formula \eqref{eq:truncatedupper}
gives
\[
 \frac1{C_m}<\frac{\widetilde{g}_{\mathbb{F}}(R_d)}{D_d}
 \leq\frac{\sum_{j=0}^d\sqrt{h_jh_{d-j}}}{D_d}.
\]
The squeeze theorem proves $D_d/\widetilde{g}_{\mathbb{F}}(R_d)\to C_m$.
Thus the supremum is $C_m$, although no individual finite
algebra attains equality in the homogeneous comparison.
\end{proof}

\begin{remark}\label{rem:singlepairing}
Regularization is essential even for three generators. For an integer
$d\geq1$, let $T_d=\mathbb{F}[x,y,z]/(x,y,z)^{2d}$. The ideal
$J=(x,y,z)^dT_d$ has square zero and codimension $\binom{d+2}{3}$.
The block argument in Proposition~\ref{prop:preliminary} gives
\[
 g_{\mathbb{F}}(T_d)\leq2\binom{d+2}{3},\qquad
 \frac{\dim_{\mathbb{F}}T_d}{g_{\mathbb{F}}(T_d)}\geq4-\frac6{d+2}.
\]
Thus a universal estimate $\dim_{\mathbb{F}}R\leq2g_{\mathbb{F}}(R)$ cannot replace
the comparison with $\widetilde{g}_{\mathbb{F}}(R)$.
\end{remark}

\section{Consequences for the generator constants}\label{sec:constants}

\subsection{Explicit coefficients}
The gamma identity $\Gamma(z+1)=z\Gamma(z)$ yields
\begin{equation}\label{eq:recurrence}
 C_{m+2}=\frac{4m}{m+1}C_m,\qquad
 C_2=\frac4\pi,\qquad C_3=2.
\end{equation}
For odd indices, $C_{2r+1}=\binom{2r}{r}$. The recurrence also gives $C_m<2^{m-1}$ for $m\geq2$: the claim holds for $m=2,3$, and the ratio to $2^{m-1}$ decreases by the factor $m/(m+1)$ on either parity subsequence. In particular, $C_m<2^m-1$.

Table~\ref{tab:coefficients} compares the coefficient in Formula \eqref{eq:oldbound} with the present upper bounds for $c_m(\mathbb{F})$. The value for $m=2$ is supplied by the commuting-pair theorem; the gamma quotient itself is $4/\pi$ in that case.

\begin{table}[htbp]
\centering
\caption{Upper coefficients for the generator constants.}\label{tab:coefficients}
\renewcommand{\arraystretch}{1.2}
\begin{tabular}{c r c r}
\toprule
$m$ & $2^m-1$ & Upper bound for $c_m(\mathbb{F})$ & Approximation\\
\midrule
2 & 3   & $1$ & $1$\\
3 & 7   & $2$ & $2$\\
4 & 15  & $32/(3\pi)$ & $3.395305$\\
5 & 31  & $6$ & $6$\\
6 & 63  & $512/(15\pi)$ & $10.864977$\\
7 & 127 & $20$ & $20$\\
8 & 255 & $4096/(35\pi)$ & $37.251351$\\
9 & 511 & $70$ & $70$\\
10 & 1023 & $131072/(315\pi)$ & $132.449249$\\
\bottomrule
\end{tabular}
\end{table}

The duplication formula for the gamma function gives
\[
 C_m=\frac{2^{m-1}}{\sqrt\pi}
 \frac{\Gamma(m/2)}{\Gamma((m+1)/2)}
 \sim\frac{2^m}{\sqrt{2\pi m}}.
\]
The improvement over the preceding coefficient $2^m-1$ is therefore
by a factor asymptotic to $\sqrt{2\pi m}$.

\subsection{Comparison with lower bounds}

\begin{proof}[Proof of Corollary~\ref{cor:growth}]
For $m\geq4$, the construction in \cite[\S3.3, equation (17)]{LLSZ}
gives $c_m(\mathbb{F})\geq L_m$ over every field. Theorem~\ref{thm:main}
gives the upper bound. The central binomial estimate implies
\[
 L_m\sim\frac{2^{m-1}}{\sqrt{2\pi m}},
\]
on both parity subsequences. Together with the preceding formula
for $C_m$, this proves $C_m/L_m\to2$.
\end{proof}

\subsection{The remaining distinction}
Theorem~\ref{thm:optimal} identifies the best coefficient in the
universal inequality $\dim_{\mathbb{F}}R\leq K_m\widetilde{g}_{\mathbb{F}}(R)$ over an
infinite field. It does not identify the best coefficient in the
faithful-module inequality $\dim_{\mathbb{F}}R\leq c_m(\mathbb{F})\dim_{\mathbb{F}}V$. The difference is already visible for two generators: $c_2(\mathbb{F})=1$, whereas the optimal regularized-rank coefficient is $C_2=4/\pi$.

For three generators, the present argument gives $c_3(\mathbb{F})\leq2$ and the finite-dimensional estimate $2n-2$, while the Gerstenhaber problem asks for $c_3(\mathbb{F})=1$. The homogeneous comparison $\dim_{\mathbb{F}}R\leq K_m\widetilde{g}_{\mathbb{F}}(R)$ alone therefore cannot yield a universal coefficient smaller than $C_m$. Further improvement requires additional information relating the algebra to its faithful modules.

    \section*{}

    \subsection*{Acknowledgement.}
         We would like to thank Prof. Kaiming Zhao for his encouragement and helpful suggestions. The second author is very grateful to Prof. Genkai Zhang for his kind hospitality and encouragement. We would also like to thank the referee for his/her detailed comments and kindness.

    \subsection*{AI disclosure.}
        During the preparation of this manuscript, the authors used LLMs to assist in drafting and checking parts of the mathematical arguments, organizing references, and improving the clarity and readability of the exposition. All mathematical statements and proofs in the final manuscript were written, reviewed, and verified by the authors, who take full responsibility for the content of this manuscript.

\end{document}